\documentclass[12pt]{article}
\usepackage{amssymb}
\makeatletter
\date{}

\newtheorem{theorem}{Theorem}[section]

\newtheorem{proposition}[theorem]{Proposition}

\newtheorem{corollary}[theorem]{Corollary}

\newtheorem{problem}[theorem]{Problem}
\newtheorem{conjecture}[theorem]{Conjecture}
\newcommand{\edim}{{\rm e}$-${\dim}}

\newcommand{\z}{{\Bbb Z}}
\newcommand{\q}{{\Bbb Q}}
\newcommand{\re}{{\Bbb R}}
\newcommand{\N}{{\Bbb N}}

\newcommand{\Int}{{\rm Int}}

\newcommand{\invlim}{{\rm invlim}}

\newcommand{\lo}{\longrightarrow}

\newcommand{\black}{{\blacksquare}}
\newcommand{\tor}{{\rm Tor}}
\newcommand{\diam}{{\rm diam}}
\newcommand{\dist}{{\rm dist}}

\begin{document}

\title{Dimension of the product and classical formulae of dimension theory}

\author{ Alexander Dranishnikov\footnote{ the first author was supported by NSF grant  DMS-0904278;}  \ and Michael Levin\footnote{the second author was supported by ISF grant 836/08}}

\maketitle
\begin{abstract} Let $f : X \lo Y$ be a map of compact metric spaces.
A classical theorem of Hurewicz asserts that
$\dim X \leq \dim Y +\dim f$
where $\dim f =\sup \{ \dim f^{-1}(y): y \in Y \}$.
The first author conjectured that {\em $\dim Y + \dim f$ in Hurewicz's theorem  can be replaced by
$\sup \{ \dim (Y \times f^{-1}(y)): y \in Y \}$}.
We  disprove this conjecture.
As a by-product of the machinery presented in the paper
we answer in negative the following problem
posed by the first author: {\em  Can for compact $X$ the Menger-Urysohn formula $\dim X \leq \dim A + \dim B +1$
be improved to $\dim X \leq \dim (A \times B) +1$ ?}

On a positive side we show that both conjectures holds true for compacta $X$ satisfying the equality $dim(X\times X)=2\dim X$.
\\\\
{\bf Keywords:} Cohomological Dimension, Bockstein Theory, Extension Theory
\bigskip
\\
{\bf Math. Subj. Class.:}  55M10 (54F45 55N45)
\end{abstract}
\begin{section}{Introduction}\label{intro}
Throughout this paper we assume that  maps are continuous and
spaces are separable metrizable. We recall that a {\em compactum}
means a compact metric space. By dimension of a space $\dim X$ we
assume the covering dimension.

Clearly, the dimension of the product of two polyhedra equals the sum of the dimension:
$\dim (K\times L)=\dim K+\dim L$.
In 1930 Pontryagin discovered that this logarithmic law does not hold for compacta~\cite{Pon}.
He constructed his famous Pontryagin surfaces $\Pi_p$
indexed by prime numbers, $\dim\Pi_p=2$, such that $\dim(\Pi_p\times\Pi_q)=3$ whenever $p\ne q$.
In the 80s the first author showed that
the dimension of the product can deviate arbitrarily from the sum of the dimension. Namely,
for any $n,m,k\in\N$ with $$\max\{n,m\}+1\le k\le n+m$$ there are compacta $X_n$ and $X_m$
of dimensions $n$ and $m$ respectively with
$\dim(X_n\times X_m)=k$~\cite{DrUsp}. We note that the inequality $\dim(X\times Y)\le\dim X+\dim Y$
always holds true.

The first author conjectured that many classical formulas
(inequalities) of dimension theory can be strengthen by replacing
the sum of the dimensions by the dimension of the product. His
believe was based on his results on the general position
properties of compacta in euclidean spaces~\cite{dranishnikov-intersection},\cite{DRS}.
Clearly, for two polyhedra $K$ and $L$ with transversal
intersection in $\re^n$ we have $\dim(K\cap L)=n-(\dim K+\dim L)$.
For compacta the corresponding formula is $\dim(X\cap
Y)=n-\dim(X\times Y)$. In particular, two compacta $X$ and $Y$ in
general position in $\re^n$ have empty intersection if and only if
$\dim(X\times Y)<n$.

The next candidate for the improvement was the following classical
theorem of Hurewicz.
\begin{theorem}  {\rm (\bf{Hurewicz Theorem})}
Let $f : X \lo Y$ be a map  of
compacta. Then
 $$\dim X \leq \dim Y +\dim f$$
where $\dim f =\sup \{ \dim f^{-1}(y) \mid y \in Y \}$.
\end{theorem}
We note that the Hurewicz theorem applied to the projection
$X\times Y\to Y$ implies the inequality $\dim(X\times Y)\le\dim
X+\dim Y$. The first author proposed the following conjecture.

\begin{conjecture}[\cite{DRS}] \label{problem1} For a map of compacta $f : X \lo Y$
$$\dim X\leq \sup \{ \dim (Y \times f^{-1}(y)) \mid y \in Y \}.$$
\end{conjecture}

Note that the Conjecture~\ref{problem1} holds true for nice maps
like locally trivial bundles. It was known that the
conjecture holds true when $X$ is standard (compactum of type I in
the sense of \cite{Ku}). We call  a compactum $X$  {\em standard}
if it has the property $\dim (X\times X)=2\dim X$. It's not easy
to come with an example of a compactum without this property. The
Pontryagin surfaces satisfy it. First example of a non-standard
compactum was constructed by Boltyanskii~\cite{Bolt}. In this paper
all non-standard compacta (compacta of type II in \cite{Ku}) will
be called {\em Boltyanskii compacta}. It is known that for all
Boltyanskii compacta $\dim(X\times X)=2\dim X-1$.

In this paper we disprove Conjecture~\ref{problem1}. We will refer
to  the maps providing counterexamples to the conjecture as {\em
exotic maps}.

Positive results towards Conjecture~\ref{problem1} can be
summarized in the following:

\begin{theorem}
\label{boltyanskii-map}
If a compactum $X$ admits an exotic map $f:X\to Y$ then
$X$ is a Boltyanskii compactum. For every exotic map $f:X\to Y$
we have $$\dim X=\sup\{ \dim (Y \times f^{-1}(y))\mid y\in Y\} +1 .$$
\end{theorem}

Another classical result
in Dimension Theory where the first author hoped to replace the sum of the
dimensions by the dimension of the product
was  the Menger-Urysohn Formula.
\begin{theorem}
{\rm ( \bf{Menger-Urysohn Formula})} Let $X=A \cup B$ be a decomposition of
a  space $X$. Then  $\dim X \leq \dim A + \dim B +1$.
\end{theorem}

\begin{problem}[ \cite{dranishnikov-intersection}]\label{problem2}
Does the inequality $\dim X \leq \dim (A\times B) +1$ hold true for an arbitrary
decomposition of compact metric space $X=A\cup B$?
\end{problem}

In this paper we
 answer  Problem \ref{problem2} in the negative and, similarly to the terminology used above,
we refer to the decompositions
providing counterexamples to Problem \ref{problem2}  as {\em exotic decompositions}.
To a certain extent exotic decompositions is a starting point of our construction
of exotic maps.

Note that in the case of non-compact $X$ a counter example to
Problem~\ref{problem2} was constructed by Jan van Mill and Roman
Pol. They proved the following.
 \begin{theorem}[\cite{mill-pol}]
 There is a $3$-dimensional subset $X \subset \re^4$ admitting
 a decomposition $X=A \cup B$ such that
 $\dim (A \times B)^n =1$ for every integer $n>0$.
 \end{theorem}
Similarly to the case of Conjecture~\ref{problem1} the following
facts were known about Problem~\ref{problem2}.

\begin{theorem}[\cite{dranishnikov-intersection})]
\label{boltyanskii-decomposition}
If a compactum $X$ admits an exotic decomposition then
$X$ is a Boltyanskii compactum. For any exotic decomposition $X=A \cup B$
of a  compactum $X$
we have $\dim X= \dim (A \times B) +2 $.
\end{theorem}
 The main results of this paper are  the following theorems.

 \begin{theorem}
 \label{general-exotic-decomposition}
 Every finite dimensional Boltyanskii compactum $X$ with $\dim X \geq 5$ admits
 an exotic decomposition.
 \end{theorem}
\begin{theorem}
\label{particular-exotic-map}
For every $n \geq 4$ there is an $n$-dimensional Boltyanskii compactum $X$
admitting an exotic map $f : X \lo Y$ to a $2$-dimensional compactum $Y$.
\end{theorem}
 Theorem \ref{particular-exotic-map}  is derived from a more general result.

\begin{theorem}
\label{initial-exotic-map}
Every   $n$-dimensional Boltyanskii compactum $X$  with $n\geq 5$ and $\dim_\q X< n-3$
admits   an exotic map $f: X \lo Y$ to an $m$-dimensional compactum $Y$ with $m=\dim_\q X +1$.
\end{theorem}

Note that no  compactum of $\dim <4$ admits an exotic map and no compactum of $\dim <5$
admits an exotic decomposition, see Section \ref{compactly-represented-spaces}.
A further development of the approach presented in the paper allows one
to partially generalize Theorem \ref{initial-exotic-map} by showing that
any finite dimensional Boltyanskii compactum $X$ with $\dim X \geq 6$
admits  an exotic map. This result is technically more complicated and
will appear elsewhere.
 It still remains open
 whether
any Boltyanski compactum of dimensions $4$ and $5$
admits an exotic map.

The paper is built as follows:  Bockstein Theory is reviewed in
Section \ref{bockstein}; Section \ref{extension-theory} is devoted to
basic facts  of  Extension Theory with applications
to Dimension Types;
in Section \ref{compactly-represented-spaces}
 we consider
the so-called compactly represented spaces, prove
Theorem \ref{general-exotic-decomposition} and
 present  short  proofs for  Theorems  \ref{boltyanskii-map}
and \ref{boltyanskii-decomposition};
and, finally,
   Theorems \ref{particular-exotic-map} and \ref{initial-exotic-map} are proved in
   Section \ref{proof-particular}.

 \end{section}

\begin{section}{Bockstein Theory}
\label{bockstein} We  recall some basic facts of Bockstein
Theory. The first detailed presentation of the theory was given in the survey \cite{Ku}. Since then it was evolved
in many papers and surveys \cite{DrUsp},\cite{DRS1},\cite{DrArxiv},\cite{dranishnikov-intersection},\cite{Sch},\cite{DyA}. Our presentation here has features
of both point of view on the subject, classical and modern.

We remind that  cohomology always means the Cech
cohomology. Let $G$ be an abelian group. The {\bf cohomological dimension} $\dim_GX$
of a space $X$ with respect to the coefficient group $G$ does not exceed $n$, $\dim _G X \leq n$ if $H^{n+1}(X,A;G)=0$ for every closed $A
\subset X$. We note that this condition implies that
$H^{n+k}(X,A;G)=0$ for all $k\ge 1$ \cite{Ku},\cite{DrArxiv}.
Thus, $\dim _G X =$ the smallest integer $n\geq 0$ satisfying
$\dim _G X \leq n$ and $\dim _G X = \infty $ if such an integer
does not exist. Clearly, $\dim_G X \leq \dim_{\z}X\le\dim X$. Note that $\dim_G
X=0$ for a non-degenerate group $G$ if and only if $\dim X=0$.

 \begin{theorem}
 \label{alexandrov}
 {\rm ({\bf Alexandroff})}
$\dim X=\dim_\z X$ if $X$ is a finite dimensional space.
\end{theorem}

Let $\mathcal P$ denote the set of all primes. The {\em  Bockstein basis} is the collection of groups
$\sigma=\{ \q, \z_p , \z_{p^\infty}, \z_{(p)} \mid p\in\mathcal P
\}$ where $\z_p =\z/p\z$ is the $p$-cyclic group,
$\z_{p^\infty}={\rm dirlim} \z_{p^k}$  is the $p$-adic circle, and
$\z_{(p)}=\{ m/n \mid n$ is not divisible by $p \}\subset\q$  is the $p$-localization of integers.
\\\\
 The Bockstein basis   of an abelian group $G$ is the collection
$\sigma(G) \subset \sigma$ determined by the rule:

$\z_{(p)} \in \sigma(G)$ if $G/\tor G$ is not divisible by $p$;

$\z_p \in \sigma(G)$ if $p$-$\tor G$ is not divisible by $p$;

$\z_{p^\infty} \in \sigma(G)$ if $p$-$\tor G\neq 0$ is  divisible by $p$;

$\q \in \sigma(G)$ if $G/\tor G\neq 0$ is  divisible by all $p$.
\\
\\
Thus  $\sigma(\z)=\{\z_{(p)}\mid p\in\mathcal P\}$.

\begin{theorem}
\label{bockstein-theorem}
{\rm ({\bf  Bockstein Theorem})} For a compactum $X$,
$$\dim_G X =\sup \{ \dim_H X : H \in \sigma(G)\}.$$
\end{theorem}

The Alexandroff and Bockstein theorems  imply that for finite dimensional compacta $X$
$$\dim
X=\max\{\dim_{\z_{(p)}}X\mid p\in\mathcal P\}.$$

We call a space $X$ {\em $p$-regular} if
 $$\dim_{\z_{(p)}} X = \dim_{\z_p} X =\dim_{\z_{p^\infty}}X=\dim_{\q}X$$
and call it {\em $p$-singular} otherwise.

The restrictions on the values of cohomological dimension of a given space
with respect to Bockstein groups usually are stated in the form of Bockstein inequalities~\cite{Ku}.
Here we state them in a form of the equality and the alternative (see \cite{DrArxiv}).
\begin{theorem}
\label{bockstein-inequalities}
I. For every $p$-singular space $X$ and  every prime $p$
$$
\dim_{\z_{(p)}}X=\max\{\dim_{\q}X,\dim_{\z_{p^\infty}}X+1\}.$$

II. (Alternative) For every $p$-singular space $X$ and  every prime $p$
either
$$ \dim_{\z_{p^\infty}}X= \dim_{\z_p} X \ \ \ or\ \ \   \dim_{\z_{p^\infty}}X= \dim_{\z_p} X-1.$$
\end{theorem}
In the first case of the alternative we call $X$ {\em $p^+$-singular} and in the second, {\em $p^-$-singular}.
Thus, the values of $\dim_{F}X$ for Bockstein fields $F\in\{\z_p,\q\}$ together with $p$-singularity types of $X$
determine the value $\dim_GX$ for all groups.

We notice that he Alexandroff theorem, the Bockstein theorem, and Theorem~\ref{bockstein-inequalities} imply the following.
\begin{corollary}
For every finite dimensional compactum $X$ there is a field $F$ such that $\dim X\le \dim_FX+1$.
\end{corollary}

A function $f:\sigma\lo \N \cup \{0, \infty\}$ is called a {\em
$p$-regular} if
$
f(\z_{(p)})=f(\z_p)=f(\z_{p^{\infty}})=f(\q)$ and it  is called {\em
$p$-singular} if
$ f(\z_{(p)})=\max\{ f(\q),f(\z_{p^{\infty}})+1\}.$
A $p$-singular function $f$ is called $p^+$-singular if
$f(\z_{p^{\infty}})=f(\z_p)$ and it is
called {\em $p^-$-singular} if
$f(\z_{p^{\infty}})=f(\z_p)-1$. A function
$D:\sigma\lo \N \cup \{0, \infty\}$ is called a {\em dimension
type} if for every prime $p$ it is either $p$-regular or
$p^{\pm}$-singular.
For every space $X$ the function $d_X:\sigma\lo \N \cup \{0,
\infty\}$ defined as $d_X(G)=\dim_GX$ is a dimension type. If $X$ is compactum
$d_X$ is called the {\em dimension type of $X$}.
We denote $\dim D=\sup \{ D(G)\mid \ G \in \sigma \}$.

\begin{theorem}
\label{dranishnikov-realization-theorem}
{\rm (\bf {Dranishnikov  Realization Theorem} \cite{DrUsp},\cite{DrPacific})}
For every   dimension type $D$ there is a compactum $X$
with $d_X=D$ and $\dim X=\dim D$.
\end{theorem}

Let $D$ be a dimension type. We will use abbreviation $D(0)=D(\q)$, $D(p)=D(\z_p)$. Additionally, if $D(p)=n\in\N$ we will write $D(p)=n^+$ if $D$ is $p^+$-regular
and $D(p)=n^-$ if it is $p^-$-regular. For $p$-regular $D$  we leave it without decoration: $D(p)=n$.
Thus, any sequence of decorated numbers $D(p)\in\N$, where $p\in\mathcal P\cup\{0\}$  define a unique dimension type.
There is a natural order on decorated numbers $$\dots<n^-<n<n^+<(n+1)^-< \dots\ .$$
Note that the inequality of dimension types $D\le D'$ as functions on $\sigma$ is equivalent to the family of inequalities
$D(p)\le D'(p)$ for the above order for all $p\in\mathcal P\cup\{0\}$.
The natural involution on decorated numbers that exchange the decorations '+' and '-' keeping the base fixed defines
an involution $\ast $ on the set of dimension types . Thus, $\ast$ takes
$p^+$-singular function $D$ to $p^-$-singular $D^*$ and vise versa.

By Alexandroff and Bockstein theorems it follows that for any compactum $X$  of $\dim X=n<\infty$ there is a prime $p$ such that $\dim X=\dim_{\z_{(p)}}X$. Then either $d_X(0)=n$ or $d_X(p)$ equals one of the following: $n$ or $n^-$ or $(n-1)^+$. Let $\mathcal P_X$ denote the set of all such primes.

The Bockstein Product Theorem \cite{Ku} gives the formulas for cohomological dimension of the product with respect to
each of the groups $G\in\sigma$ which are huge for some of $G$. Here we state it in an alternative way (see \cite{DrArxiv},\cite{Sch},\cite{DyA}).
\begin{theorem}
\label{bockstein-product-theorem}
{\rm(\bf {Bockstein Product Theorem})}
For every field $F$ and any two compacta, $$\dim_{F}(X \times Y)=\dim_{F} X + \dim_{F} Y.$$
For every prime $p$ the type of $p$-singularity is preserved by multiplication by a $p$-regular compactum, and
the following rule is applied in the remaining cases:

$p^+$-singular $\times$ $p^+$-singular = $p^+$-singular;

$p^-$-singular $\times$ $p^{\pm}$-singular= $p^-$-singular.
\end{theorem}
The product formula implies that an $n$-dimensional compactum $X$ is a Boltyanskii compactum if and only if
$d_X(0)<n$ and $d_X(p)=(n-1)^+$ for all $p\in\mathcal P_X$ . For every $n\ge 2$ we denote by $B_n,$ the "maximal" dimension type of Boltyanskii compacta of dimension $n$.
Thus, $B_n(p)=(n-1)^+$ for all $p\in\mathcal P$ and $B_n(\q)=n-1$. This implies that $B_n(\z_{(p)})=n$ for every
prime $p$ and $B_n(G)=n-1$ for all other groups in $\sigma$.
\begin{corollary}
\label{boltyanskii-type}
{\rm ({\bf \cite{dranishnikov-intersection}})}
For an $n$-dimensional compactum $X$ the following are equivalent:
\begin{itemize}
\item  $X$ is a Boltyanskii compactum;
\item $\dim X > \dim_F X$ for every field $F$;
\item $d_X(G) \leq B_n(G)$ for all $G\in\sigma$.
\end{itemize}

A finite dimensional compactum $X$ is   standard   if and only if  there  is
a field $F\in \sigma$ such that $\dim X= \dim_F X$.
\end{corollary}

Let $D_1$ and $D_2$ be dimension types.
 The dimension
type  $D_1 \boxplus D_1$ is defined   by
 the formulas of the Bockstein Product Theorem:
 $(D_1 \boxplus D_2) (G) = \dim_G(X \times Y)$ with $\dim_G X$ and $\dim_G Y$
 being replaced by  $D_1(G)$ and $D_2(G)$ respectively for $G \in \sigma$ (see~\cite{DRS1}).
 Thus we have that
 $d_{X \times Y}=d_X \boxplus d_Y$ for compacta $X$ and $Y$. If $D_1(p)=n^{\epsilon_1}$ and $D_2(p)=m^{\epsilon_2}$
where $\epsilon_i$ is a decoration, i.e., '+' or '-' or empty, then $$(D_1\boxplus D_2)(p)=(n+m)^{\epsilon_1\otimes\epsilon_2}$$
with the product of the signs $\epsilon_1\otimes\epsilon_2$ defined by
the Bockstein Product Theorem rule:
$$\epsilon\otimes empty=\epsilon,\ \  \ \ \epsilon\otimes \epsilon=\epsilon, \ \ \epsilon=\pm,\ \ \  and\ \ \
+\otimes -=-.$$

 By   $D_1 +D_2$ and $D_1 \leq D_2$ we mean
 the ordinary   sum and order relation  when   $D_1$ and $D_2$  are
 considered as just functions.
 Note  that   $D_1 +D_2$ is not always
 a  dimension type but it is a dimension type, provided one of the summands is $p$-regular for all $p$.
 By $0$ and $1$ we denote the dimension types which send every
 $G\in \sigma$ to $0$ and $1$ respectively. Recall that $d_X=0$
 if and only if $\dim X = 0$ and
   $$d_{X\times [0,1]}=d_X +1.$$

The following inequality is an  easy observation.
\begin{proposition}\label{oplus} For any dimension types $D_1$ and $D_2$,
$$ D_1\boxplus D_2\le (D_1^*\boxplus D_2^*)^*.$$
\end{proposition}
{\bf Proof.} Clearly, we have the equality $ (D_1\boxplus D_2)(F)= (D_1^*\boxplus D_2^*)^*(F)$
for the fields. Thus, it suffices the check the inequality for
the decorations. If $D_1^*\boxplus D_2^*$ is $p^-$-singular,
then the right hand part will have the decoration '+' and the inequality holds. If $D_1\boxplus D_2$ is $p^-$ singular,
clearly the inequality holds. In the remaining case both $D_1$ and $D_2$ are $p$-regular and therefore, we have the equality.
$\black$

\end{section}
\begin{section}{Extension Theory}
\label{extension-theory}
Cohomological Dimension is characterized by the following basic property:
$\dim_G X \leq n$ if and only  for every closed
$A \subset X$ and a map $f : A \lo K(G,n)$,
$f$ continuously extends over $X$ where $K(G,n)$ is the Eilenberg-MacLane complex
of type $(G,n)$
(we assume that $K(G,0)=G$ with discrete topology and $K(G, \infty)$ is a singleton).
This extension characterization of Cohomological
Dimension gives a rise to  Extension  Theory (more general than
Cohomological Dimension Theory)
 and  the notion of Extension Dimension.
 The {\em extension dimension} of a space $X$ is said
to be dominated by a CW-complex $K$, written $\edim X \leq K$, if
every map $f : A \lo K$ from a closed subset $A$ of $X$
continuously extends over $X$. Thus $\dim_G X \leq n$ is equivalent
to $\edim X \leq K(G,n)$ and $\dim X \leq n$ is equivalent to $\edim X \leq S^n$.
For a dimension type  $D$ we denote
 $K(D)=\vee_{G \in \sigma} K(G, { D}(G))$.
 Then  $d_X \leq D$ if and only
 $\edim X \leq K(D)$.
\\\\
Extension Dimension has many properties similar to Covering Dimension.
For example:
if $\edim X \leq K$ then $\edim A \leq K$ for every
$A \subset X$ and if $X=\cup F_i$ is a countable union of
closed subsets of $X$ such that $\edim F_i\leq K $ for every $i$
then $\edim X \leq K$.
Let us list  a few more basic results of Extension Theory.

\begin{theorem}
\label{olszewski-completion-theorem}
{\rm (\bf {Olszewski Completion  Theorem} \cite{Ol})}
Let $K$ be a countable CW-complex and $\edim X \leq K$.
Then there is a completion of $X$ dominated by $K$.
\end{theorem}
\begin{corollary}
For every  separable metric space $X$ there is a completion $X'$ such that for all $G\in\sigma$, $$\dim_GX'=\dim_GX.$$
\end{corollary}
We note that for finite dimensional $X$ this corollary follows from the theory of test spaces~\cite{Ko},
\cite{Ku}
and the well-known fact that for every compactum $C$ there is a completion $X'$ of $X$ with
$\dim(X'\times C)=\dim(X\times C)$ (see for example Proposition 6.2 in\cite{dranishnikov-intersection}).

\begin{theorem}
\label{dranishnikov-extension-theorem}
{\rm (\bf {Dranishnikov Extension Theorem} \cite{DrExt},\cite{Dy})}
Let $K$ be a  CW-complex and
  $X$ a space. Then

(i)
$\dim_{H_n(K)} X \leq n$ for every $n\geq 0$
 if
$\edim X \leq K$;

(ii) $\edim X \leq K$ if $K$ is simply connected, $X$ is finite dimensional
and
$ \dim_{H_n(K)} X \leq n$ for every $n\geq 0$
\end{theorem}
We remind that $H_*(K)$ denotes the reduced  homology.
\\\\
Let $K$ be a CW-complex.  For $G \in \sigma$ denote
$n_G(K) =\min \{ n : G \in \sigma(H_n(K)) \}$ or
$n_G(K) =\infty$ if the set defining $n_G(K)$  is empty.
If $X$ is a compactum and $\edim X \leq K$ then,
by the Dranishnikov Extension Theorem and the Bockstein Theorem, we have that
$\dim_{G} X \leq n_G(K)$ for every $ G\in \sigma$.

\begin{theorem}
\label{dranishnikov-union-theorem}
{\rm (\bf {Dydak Union  Theorem} \cite{Dy})}
Let $K$ and $L$ be CW-complex and $X=A \cup B$
a decomposition of a space $X$ such that $\edim A \leq K$
and $\edim B \leq L$. Then $\edim X \leq K*L$.
\end{theorem}
We recall that $K*L=\Sigma (K\wedge L$).

\begin{theorem}
\label{dranishnikov-decomposition-theorem}
{\rm (\bf {Dranishnikov Decomposition  Theorem} \cite{DrPacific})}
Let $K$ and $L$ be countable CW-complexes and $X$ a compactum such that
$\edim X \leq K*L$. Then there is  a decomposition  $X=A \cup B$
 of  $X$ such that $\edim A \leq K$
and $\edim B \leq L$.
\end{theorem}
Let $D_1$ and $D_2$ be   dimension types such that
at least one of them is different from $0$ and $X=A \cup B$ a decomposition
of a compactum $X$ such that $d_A\leq D_1$ and $d_B\leq D_1$.
By the Dydak Union Theorem, $d_X \leq K(D_1)*K(D_2)$. Then
 $\dim_G X \leq n_G(K(D_1) * K(D_2))=
 n_G(\Sigma(K(D_1) \wedge K(D_2)))=n_G(K(D_1) \wedge K(D_2)) +1$, $G\in \sigma$.
 Thus
 one can estimate the dimension type of $X$ by computing
 the numbers $n_G(K(D_1) \wedge  K(D_2))$, $G\in \sigma$.
 This computation was done by Dranishnikov \cite{dranishnikov-intersection}.
 We denote by
       $D_1 \oplus D_2$  the biggest dimension type  such that
  $(D_1 \oplus D_2)(G) \leq n_G(K(D_1) \wedge K(D_1)), G \in \sigma$ and
   set   $D_1\oplus D_2=0$ for $D_1=D_2=0$.
  The following can be easily derived from Dranishnikov's computation~\cite{dranishnikov-intersection}:

 \begin{theorem}
 \label{smash-product-extension-type}
Let $D_1$ and $D_2$ be  dimension types. Then
 $$D_1 \oplus D_2=(D_1^*\boxplus D_2^*)^*. $$
\end{theorem}
Thus if  $X=A \cup B$ is a decomposition of a compactum $X$ with
$d_A \leq D_1$ and $d_B \leq D_2$ then
$d_X \leq D_1\oplus D_2 +1$.
\\\\
Now assume that $X$ is a finite dimensional compactum
and  $D_1$ and $D_2$ are  dimension types such that
$d_X \leq D_1 \oplus D_2 +1$. If  $D_1=D_2=0$ then $\dim X\leq 1$
and for a decomposition  $X=A \cup B$ into $0$-dimensional subsets
we obviously have $d_A\leq D_1$ and $d_B  \leq  D_2$.  If at least one
of $D_1$ and $D_2$ is different from $0$ then  $K(D_1)*K(D_2)$
is simply connected. Then, by  the Bockstein Theorem
 and  the Dranishnikov Extension Theorem,
 $\edim X \leq K(D_1)*K(D_2)$ and, by the Dranishnikov Decomposition Theorem,
 there is  a decomposition $X=A \cup B$ of $X$ with $\edim A \leq K(D_1)$
 and $\edim B \leq K(D_2)$ and, hence, $d_A \leq D_1$ and $d_B\leq D_2$.
 Thus we can summarize

 \begin{corollary}
 \label {union-decomposition-for-dimension-types}
 Let $X$ be a compactum and $D_1$ and $D_2$  dimension types:

 (i) if $X=A \cup B$ is  a decomposition with $d_A \leq D_1$ and
 $d_B \leq D_2$ then $d_X\leq D_1 \oplus D_2 +1$;

 (ii) if $X$ is finite dimensional and  $d_X \leq D_1 \oplus D_2 +1$
 then there is a decomposition  $X=A \cup B$ such that
 $d_A \leq D_1$ and   $d_B \leq D_2$.

\end{corollary}
Note that
$$(D_1 \boxplus D_2)(F) =(D_1 \oplus D_2)(F) = D_1(F) + D_2(F)$$
for any dimension types $D_1$, $D_2$ and any field $F \in \sigma$.
\begin{proposition}\label{compare}
Let the dimension types $D_1,D_2, D'_1$ and $D'_2$ satisfy
$D_1 \leq D'_1$ and $D_2 \leq D'_2$. Then
$D_1 \boxplus D_2 \leq D'_1 \boxplus D'_2$ and
$D_1 \oplus D_2 \leq D'_1 \oplus D'_2$
\end{proposition}
{\bf Proof.}
The first inequality is standard and it easy follows from the definitions. The second inequality follows
from Theorem~\ref{smash-product-extension-type}.
$\black$

It turns out that the operation $\oplus$ nicely fits in the translation of some
mapping theorems by Levin and Lewis   \cite{levin-lewis} to
the language of dimension types.

\begin{theorem}
\label{levin-lewis-domain}
{\rm ({\bf Levin-Lewis \cite{levin-lewis}})}
Let $f: X \lo Y$ be a map of compacta and let
$K$ and $L$ be CW-complexes such that
$\edim f  \leq K$ and $\edim Y \leq L$.
Then $X \times [0,1]$ decomposes into
$X \times [0,1]=A \cup B$ such that $\edim A \leq K$ and
$\edim B \leq L$.
\end{theorem}

\begin{theorem}
\label{levin-lewis-fibers}
{\rm ({\bf Levin-Lewis \cite{levin-lewis}})}
Let $f:  X \lo Y$ be a map of compacta and,
$K$ a countable  CW-complexes such that
$\edim f  \leq \Sigma K$ and $Y$ is finite dimensional. Then

(i)
 there is a $\sigma$-compact  set $A \subset X$
such that $\edim A \leq K$ and $ \dim f|_{X \setminus A} \leq 0$;

(ii) there is a map $g : X \lo [0,1]$ such that for
the map $(f,g) : X \lo Y \times [0,1]$ we have
$\edim(f,g) \leq K$.
\end{theorem}
Let $f : X \lo Y$ be  a map.
For a group $G$ we denote $\dim_G f=\sup \{ \dim_G f^{-1}(y): y \in Y\}$ and
for  a CW-complex $K$
we say that $\edim f \leq K$ if
 $\edim f^{-1}(y) \leq K$ for every $y \in Y$. Similarly, for a dimension type $D$
we say that $d_f \leq D$ if  $d_{f^{-1}(y)} \leq D$ for every $y \in Y$.
Theorems  \ref{levin-lewis-domain} and
\ref{levin-lewis-fibers}
can be translated  to dimension types as follows.
\begin{corollary}
\label{dimension-type-domain}
Let $f: X \lo Y$ be a map of compacta and let
$D_1$ and $D_2$ be  dimension types  such that
$d_f \leq D_1$ and $d_Y \leq D_2$.
Then $d_X \leq D_1 \oplus D_2$. Moreover, if $F \in \sigma$ is a field
then $\dim_F X \leq \dim_F f +\dim_F Y$.
\end{corollary}
{\bf Proof}. Apply Theorem \ref{levin-lewis-domain}
for $K=K(D_1)$  and $L=K(D_2)$ to get
a decomposition $X \times [0,1]= A \cup B$
with $\edim A \leq K(D_1)$ and $\edim B \leq K(D_2)$.
Then $d_A \leq D_1$, $d_B \leq D_2$ and,
by  \ref{union-decomposition-for-dimension-types},
$d_X +1=d_{X \times [0,1]}\leq D_1 \oplus D_2 +1$ and
hence $d_X \leq D_1 \oplus D_2$.

 Now let $F\in \sigma $ be a field, $n=\dim_F f$, $m =\dim_F Y$ and let
 $K=K(F,n)$ and $L=K(F,m)$. Then, by the reasoning we just used,
there is a decomposition  $X \times [0,1]= A \cup B$ such that
 $\edim A \leq K$ and $\edim B \leq L$ and, by Corollary~\ref{union-decomposition-for-dimension-types},
 $d_X +1\leq d_A \oplus d_B +1$. Hence
 $\dim_F X =d_X(F)\leq (d_A \oplus d_B)(F)=
  d_A (F) + d_B (F)=\dim_F A +\dim_F B \leq n+m=
  \dim_F f + \dim_F Y$.
 $\black$

\begin{corollary}
\label{dimension-type-fibers}
Let $f:  X \lo Y$ be a map of finite dimensional compacta and
$D$ a  dimension type  such that
$d_f  \leq D+1$. Then

(i) there is a $\sigma$-compact  set $A \subset X$ such that $ d_A
\leq  D$ and $\dim (f|_{X \setminus A}) \leq 0$;

(ii) there is a map $g : X \lo [0,1]$ such that
for the map $(f,g) : X \lo Y\times [0,1]$ we have
$d_{(f,g)} \leq D$.
\end{corollary}
{\bf Proof.} By Corollary \ref{union-decomposition-for-dimension-types}
we have that each fiber $f^{-1}(y)$ decomposes into
$f^{-1}(y)=\Omega_1 \cup \Omega_2 $ with  $d_{\Omega_1} \leq 0$  and $d_{\Omega_2} \leq D$.
Then $\edim \Omega_1 \leq S^0$, $\edim \Omega_2 \leq K(D)$ and,
by the Dydak Union Theorem, $\edim f^{-1}(y) \leq S^0 * K(D) =\Sigma K(D).$

Thus,  $\edim f \leq \Sigma  K(D)$ and  the corollary follows from
Theorem \ref{levin-lewis-fibers}. $\black$
\\\\
We end this section with the following observation.
\begin{proposition}
\label{rational-dimension}
Let $X$ be a finite dimensional compactum and $n>0$.
Then  $\dim_\q  X \leq n$ if and only if for every closed subset
$A$ of $X$ and every map $f : A \lo S^n$ there is a map
$g : S^n \lo S^n$ of non-zero degree such that
$g \circ f : X \lo S^n$ continuously extends over $X$.
\end{proposition}
{\bf Proof.} Let $M(\q,n)$ be a Moore space of type $(\q,n)$.
Represent $M(\q, n)$ as  the telescope of a sequence  of maps
$\phi_i : S^n \lo S^n$ such that $\deg \phi_i=i, i> 0$.
Note that $M(\q,1)=K(\q,1)$. By
the Dranishnikov Extension Theorem $\edim X \leq M(\q, n)$ is equivalent to
$\dim_\q X \leq n$ for $n \geq 2$.
Thus  $\edim X \leq M(\q, n)$ is equivalent to $\dim_\q X \leq n$ for every $n>0$.

Assume that $\dim_\q X \leq n$.
Consider $f$ as a map to the first sphere of $M(\q, n)$ and continuously
extend $f$ to $f ' : X \lo M(\q,n)$. Then $f'(X)$ is contained in a finite
subtelescope $M'$ of $M(\q, n)$. Let $r : M' \lo S^n$ be  the natural retraction
to the last sphere of $M'$. Then $g$ can be taken as $r$ restricted
to the first sphere of $M'$.

Now  we will show the other direction of the proposition.
Take  a map  $\psi : A \lo M(\q,n)$   from a closed subset
$A$ of $X$. Then $\psi(A)$ is contained in a finite subtelescope
of $M(\q, n)$. Assume that that this subtelescope ends
at the $i$-th sphere of $M(\q, n)$.
Then
$\psi$ can be homotoped to a map $f: A \lo S^n$ to
the $i$-th sphere of  $M(\q, n)$. Let $g : S^n \lo S^n$ be
a map of degree $d>0$ such that $g \circ f $ extends over $X$. Consider
the subtelescope $M'$ of $M(\q, n)$ starting at the $i$-th sphere
and ending at the $(i+d)$-sphere of  $M(\q, n)$ and
let $r : M ' \lo S^n $ be the natural retraction of $M'$ to the last sphere of
$M'$.
 Then the degree of $r$ restricted to the $i$-the sphere  of $M(\q, n)$
 is divisible by $d$
 and hence $r\circ f$ factors up to homotopy through $g \circ f$. Since
 $r\circ f$ and $f$ are homotopic as maps to $M'$ we get that $f$ extends
 over $X$ as a map to $M'$ and therefore $\psi$ extends as well.
 $\black$

\end{section}

\begin{section}{Proofs of Theorems \ref{boltyanskii-map}, \ref{boltyanskii-decomposition}
and \ref{general-exotic-decomposition}}
\label{compactly-represented-spaces}
A space $X$ is called {\em compactly represented } if
for every $G \in \sigma \cup \{ \z \}$ there is a compactum
$C \subset X$ such that $\dim_G C=\dim_G X$.
 We say that
a  space $X$ is {\em compactly represented by a subset} $A \subset X$ if
$X$ is  compactly represented and the compacta $C$ witnessing that
can be chosen to be subsets of $A$. Note that
any $\sigma$-compact set is compactly represented.
 We say that a space $X$ is {\em dimensionally dominated by a space} $Y$
if $\dim_G X \leq \dim_G Y$ for every $G \in \sigma \cup \{ \z \}$.
It  follows from
the Olszewski Completion Theorem
that for a $\sigma$-compact subset  $A $ of a compactum $X$
and a space $Y$ there is a $G_\delta$-subset  $A'\subset X$
such that $A \subset A'$, $A'$ is compactly
represented by $A$ and $A'\times Y$ is dimensionally dominated by $A\times Y$.
Moreover, if $Y$ is also $\sigma$-compact then we may assume that
$A'\times Y$ is compactly represented by $A\times Y$.
Note that
$\dim_\z X =\sup \{ \dim_G  X : G \in \sigma \}$ if $X$ is compactly represented
 and $d_{X\times Y} =d_X \boxplus d_Y$ if $X$, $Y$ and
$X \times Y$ are compactly represented.
 \\\\
We say that  a decomposition $X=A \cup B$ of a space $X$ is a {\em compactly represented decomposition} if
 $A$, $B$ and $A \times B$ are compactly
 represented and
  we say that
 a decomposition $X=A' \cup B'$ is   {\em dimensionally dominated
 by  a decomposition} $X=A\cup B$ if   $\dim_G A \leq \dim_G A'$, $\dim_G B \leq \dim_G B'$
 and $\dim_G (A\times B) \leq \dim_G(A'\times B')$ for every $G \in \sigma \cup \{ \z \}$.
 \\\\
The following proposition can be  easily derived from
the proof of Proposition 6.3 of \cite{dranishnikov-intersection}.

\begin{proposition}
\label{dranishnikov-intersection-domination}
Let $X$ be a  compactum, and $X=A \cup B$ a decomposition. Then
there is a decomposition $X=A' \cup B'$ such that $A'$ is $\sigma$-compact,
$B'=X\setminus A'$
and the decomposition $X=A'\cup B'$    is dimensionally dominated
by the decomposition $X=A \cup B$.

\end{proposition}
 We need a stronger version of Proposition \ref{dranishnikov-intersection-domination}.
 \begin{proposition}
 \label{compactly-represented-domination}
 Let $X$ be a compactum. For any decomposition $X =A \cup B$
 of $X$ there is a compactly represented  decomposition $X=A' \cup B'$
 such that $A'$ is $\sigma$-compact, $B'=X\setminus A'$ and
 the decomposition $X=A'\cup B'$
 is dimensionally dominated by the decomposition $X=A \cup B$.

 \end{proposition}
 {\bf Proof.} By Proposition \ref{dranishnikov-intersection-domination},
 we can assume that $B$ is $\sigma$-compact and $A=X \setminus B$.
  Let $B_1$ be a $G_\delta$-subset of $X$ such that
  $B \subset B_1$, $B_1$ is compactly represented by $B$
  and $A\times B_1$ is dimensionally dominated by $A \times B$.
  Set $A_1=X \setminus B_1$. Then there is
  a $G_\delta$-subset $B_2$ of $X$ such that $B  \subset B_2\subset B_1$,
  $B_2$ is compactly represented by $B$ and $A_1 \times B_2$ is compactly
  represented by $A_1\times B$.
     Proceed by induction and construct for every
   $i$  a $G_\delta$-set  $B_i $ and  a $\sigma$-compact set $A_i=X \setminus B_i$
   such that \\

  (i) $B \subset B_{i+1} \subset B_i$;

(ii) $B_i$ is compactly represented by $B$;

(iii) $A_i \times B_{i+1}$ is compactly represented by $A_i\times B$.
\\
\\
Then for $B'=\cap B_i$  and $A'=\cup A_i$ we have
that $X=A' \cup B'$, $A' \subset A$
$B \subset B'$, $A'$ is $\sigma$-compact, $B'$ is $G_\delta$,
$B'$ is compactly represented by $B$, $A'\times B'$ is
compactly represented by $A' \times B$. Recall that
$A_i \times B' \subset A \times B_1$  and $A \times B_1$ is dimensionally
dominated by $A \times B$. Thus
$A' \times B'$  is dimensionally dominated by $A\times B$
and the proposition follows.
$\black$
\\

{\bf Proof of Theorem \ref{boltyanskii-map}.}  Let $f : X \lo Y$ be an exotic map of
compacta. Then for every field $F \in \sigma$ and every $y \in Y$ we have
$$\dim X >\dim (f^{-1}(y) \times Y) \geq \dim_F( f^{-1}(y) \times Y) =\dim_F f^{-1}(y) + \dim_F Y.$$
Hence $$\dim X >  \sup_{y\in Y}\{\dim (f^{-1}(y)\times Y)\}\geq
\sup_{y\in Y}\{\dim_F (f^{-1}(y)\times Y)\} +\dim_F Y=\dim_F f +\dim Y.$$
Then,
by Corollary \ref{dimension-type-fibers},  $ \dim_F f +\dim Y \geq \dim_F X$ and
hence $$\dim X >  \sup\{\dim (f^{-1}(y)\times Y): y \in Y \}\geq \dim_F X$$
for every field $F \in \sigma$.
Thus, by Corollary \ref{boltyanskii-type},
we conclude that $X$ is a Boltyanskii compactum and
$\dim X =\sup\{\dim (f^{-1}(y)\times Y)\mid\ y \in Y \}+1$. $\black$
\\
\\
{\bf Proof of Theorem \ref{boltyanskii-decomposition}.}
let  $X=A \cup B$ be an exotic decomposition of a compactum $X$.
By Proposition \ref{compactly-represented-domination} we may assume that
 $X=A \cup  B$ is a compactly represented decomposition. Then for every field
 $F \in \sigma$ we have $$\dim_F X \leq \dim_F A + \dim_F B +1 =
 \dim_F (A \times B) +1 \leq \dim (A \times B) +1 \leq \dim X - 1.$$
 Thus $\dim X \geq \dim_F X +1$  for every field $F \in \sigma$.
 Then, by Corollary \ref{boltyanskii-type}, $X$ is a Boltyanskii compactum
 and
  there is a field $F$ such that $\dim_F X +1=\dim X$. Hence
 $\dim X = \dim (A \times B) +2$ and the theorem follows. $\black$
\\
\\
{\bf Proof of Theorem \ref{general-exotic-decomposition}.}
Let $X$ be an $n$-dimensional Boltyanskii compactum with $n \geq 5$.
Define the dimension types $D_1$ and $D_2$ by
$D_1(p)=2^-$, $D_1(\q)=1$ and
$D_2(p)=(n-4)^+$, $D_2(\q)=n-3$
for all primes $p$.
Then $$(D_1 \oplus D_2)(p) =(2^+\boxplus (n-4)^-)^*=((n-2)^-)^*=(n-2)^+$$ for all $p$
and $(D_1\oplus D_2)(\q)=D_1(\q)+D_2(\q)=n-2$. Thus, $D_1\oplus D_2=B_{n-1}$
where $B_{n-1}$ is the maximal Boltyanskii dimension type of dimension $n-1$.
Since $(D_1\boxplus D_2)(p)=(n-2)^-$ and $(D_1\boxplus D_2)(\q)=n-2$, we obtain
that $(D_1\boxplus D_2)(\z_{(p)})=n-2$ and hence,
$\dim (D_1 \boxplus D_2)   \leq n-2 $.
  By Corollary \ref{boltyanskii-type},
 $d_X \leq B_n =B_{n-1} +1 =D_1\oplus D_2 +1$ and,
 by Corollary \ref{union-decomposition-for-dimension-types},
 there is a decomposition $X=A \cup B$ such that
 $d_A \leq D_1$ and $d_B \leq D_2$.
 By Proposition \ref{compactly-represented-domination} we can assume
 that  $X=A \cup B$ is a compactly represented decomposition.
 Then $d_{A \times B} \leq D_1 \boxplus D_2$ and
  $$\dim (A \times B)=\dim_\z( A \times B) =\dim d_{A\times B}\le\dim( D_1 \boxplus D_2)   \leq n-2.$$
 Thus
 $X=A \cup B$ is an exotic decomposition
 and the theorem follows. $\black$
 \\
 \\
 Note that for compacta  $X$ and $Y$ with $\dim Y \geq 1$ we always
 have $\dim (X\times Y)\geq \dim X +1$. This property immediately
 implies that no compactum of dimension$\leq 3$
 admits an exotic map. Together with
 Proposition \ref{compactly-represented-domination} this property
 also implies that no compactum of dimension$\leq  4$ admits
 an exotic decomposition.
\end{section}

\begin{section}{Proofs of Theorems \ref{particular-exotic-map} and \ref{initial-exotic-map}}
\label{proof-particular}

For  $n \geq 5$ and $n-3\geq m\geq 2$
consider the dimension types
$D, D_1$ and $D_2$ defined by  $D(\q)=D_1(\q)=m-1, D_2(\q)=n-m-1$ and for every  $p\in\mathcal P$,
$$D(p)=(n-1)^+,\ \ \ \ \ D_1(p)=m^-,\ \ \ \ \ D_2(p)=(n-m-2)^+.$$
Note that  $(D_1\oplus D_2)(p)=(n-2)^+$ and $(D_1\boxplus D_2)(p)=(n-2)^-$. Hence,
$D\leq D_1 \oplus D_2 +1$ and $\dim (D_1\boxplus (D_2+1)) =n-1$.

Note that for an $n$-dimensional Boltyanskii compactum with $n \geq 5$ and $m=\dim_\q X +1 \leq n-3$
we have $d_X \leq D$. Then Theorem \ref{initial-exotic-map} immediately follows
from  the following proposition.
\begin{proposition}
\label{first-exotic-map}
Every $n$-dimensional compactum $X$ with
 $d_X \leq D$  admits a map
 $f : X\lo Y$  to an $m$-dimensional compactum $Y$ such that
$d_Y  \leq D_1$ and $d_f \leq D_2+1$. Thus for every
$y \in Y$,
$$\dim (f^{-1}(y)\times Y) \leq \dim (D_1\boxplus (D_2 +1)) =n-1$$ and
hence $f$ is an exotic map.
\end{proposition}
All the cases of Theorem \ref{particular-exotic-map},  except $n=4$, are covered
by Theorem \ref{initial-exotic-map} for $m=2$.
Let us  show that the missing case  $n=4$
also follows from Proposition \ref{first-exotic-map}.
\\
{\bf Proof of Theorem \ref{particular-exotic-map} (the missing  case).}
Consider  the map $f : X \lo Y$  constructed in
Proposition \ref{first-exotic-map} for $n=5$ and $m=2$.
By Theorem \ref{dimension-type-fibers},
there is a map $g : X \lo [0,1]$ such that the map $(f,g) : X \lo Y \times [0,1]$
is of dimension type $d_{(f,g)} \leq D_2$. By the Hurewicz Theorem there is $t \in [0,1]$ such that
$X'=g^{-1}(t)$ is of $\dim \geq 4$. Let $f': X' \lo Y$ be the map $(f,g)|_{X'}$ followed
by the projection from $Y \times [0,1]$ to $[0,1]$ and
let $Y'=f'(X')$. Then $d_{f'} \leq  D_2$
and $d_{Y'}\leq D_1$.
Thus $\dim f' \leq 2$ and $\dim Y' \leq 2$ and since $\dim X' \geq 4$ we get,
by the Hurewicz Theorem, that $\dim X'=4$ and $\dim f' =\dim Y'=2$.
 Note that
$\dim(D_1 \boxplus D_2) =3$ and hence   $f' : X'\lo Y'$  is an exotic map
we are looking for.
 $\black$
 \\\\
 In the proof of Proposition \ref{first-exotic-map}
 we will use the following.
\begin{proposition}
\label{known} Let $X$ be a compactum, $M$ an $m$-dimensional
manifold possibly with boundary, $A$ a $\sigma$-compact subset of
$X$ with $\dim A \leq m$, $F$ a closed subset of $X$ and $f : X
\lo M$ a map which is $0$-dimensional on $A\cap F$. Then $f$ can
be arbitrarily closely approximated by a map $f' : X \lo M$ such
that $f'$ is $0$-dimensional on $A$ and $f'$ coincides with $f$ on $F$.
\end{proposition}
{\bf Proof.} This a typical application of the Baire Category
Theorem. The following fact is well known:

{\em (1) For every $m$-dimensional compactum $A$ the set $G$ of
$0$-dimensional maps $f:A\to\re^m$ is dense $G_{\delta}$ in the
space of all continuous maps $C(A,\re^m)$ given the uniform
convergence topology.}

For the proof we  present $G$ as the intersection of sets $W_n$ of
maps $f:A\to \re^m$ such that $diam\ C<1/n$ for all components $C$
of the preimage $f^{-1}(x)$ for all $x$. It is easy to see that
each $W_n$ is open. One way to show that $W_n$ is dense is first
to approximate a given map $f:A\to\re^m$ by a composition $f'\circ
q$ where $q:A\to K^m$ is an $(1/n)$-map to an $m$-dimensional
simplicial complex and then to approximate $f'$ by a
$0$-dimensional map $g:K^m\to\re^m$. The latter can be obtained by
a proper perturbation of all vertices of sufficiently small
subdivision of $K^m$ in $\re^m$ and by taking the corresponding
piece-wise linear map $g$.

We note that $\re^m$ can be replaced by the half space $\re_+^m$
in this proof.

The above statement can be generalized to the following:

{\em (2) For every compact metric pair $(X,A)$ with
$m$-dimensional $A$ the set $G$ of maps $f:(X,A)\to (M,B)$  with
$0$-dimensional restriction $f|_A$ is dense $G_{\delta}$ in the
space of all continuous maps of pairs $C((X,A),(M,B))$  where $M$
is compactum and $B\subset M$ is an open set homeomorphic to
$\re^m$ or to $\re^m_+$.}

Note that $C((X,A),(M,B))$ is open in $C(X,A)$ and hence is
complete. Then the above argument works for this statement as
well.

As a corollary we obtain the following:

{\em (3) For every compact metric pair $(X,A)$ with
$m$-dimensional $A$ and an $m$-dimensional manifold with boundary
$M$ every continuous map $f:X\to M$ can be approximated by maps
$0$-dimensional on $A$.}

To derive it from the above we consider a finite cover $B_1,\dots
B_k$ of $f(A)$ by open sets homeomorphic to $\re^m$ or $\re^m_+$
and a partition $A=A_1\cup\dots\cup A_k$ into closed subsets such
that $f(A_i)\subset B_i$. Let
$$W=W(\{A_i\},\{B_i\})=\{g:X\to M\mid g(A_i)\subset B_i,
i=1,\dots,k\}$$ be a corresponding neighborhood of $f$ in the
compact-open topology. Then the set $W\cap(\cap_iG_i)$ is dense
$G_{\delta}$ in $W$ where $G_i$ the set of maps $g:(X,A_i)\to (M,B_i)$ which
are $0$-dimensional on $A_i$.

We note that the compactness of $A$ in this statement can be
replaced by $\sigma$-compactness.

\

Finally, to obtain the statement of the proposition we consider a
compact subset $F\subset X$ with a fixed map $f_0:F\to M$ which is
$0$-dimensional on $F\cap A$. Let $A'=A\setminus F$. Note that
$A'=\cup A_i$ is the countable union of compact sets $A_i$. Now we
prove the statement (3) for $A'$ in the complete metric space
$C(X,M;F,f_0)=\{f:X\to M\mid f|_F=f_0\}$ using the same proof. By
the countable  union theorem any map $f\in C(X,M;F,f_0)$ which is
$0$-dimensional on $A'$ is also $0$-dimensional on $A$.
$\black$\\\\

{\bf Proof of Proposition \ref{first-exotic-map}.}
 Since $D\leq D_1 \oplus D_2 +1$, by Corollary \ref{union-decomposition-for-dimension-types}
 there is a decomposition $X=A \cup B$ such that $d_A \leq D_1$ and  $d_B \leq D_2$.
 By the corollary of Olszewski Completion Theorem we may assume that $B$ is $G_\delta$. Thus replacing
 $A$ by $X \setminus B$  we assume that $A$ is $\sigma$-compact.
 Represent $A=\cup A_i$ as a countable union compact subsets $A\subset X$
 such that $A_i \subset A_{i+1}, i=1,2,\dots$. Note that $\dim A \leq \dim D_1=m$.
 \\
 \\
 We will construct for each $i$
an  $m$-dimensional simplicial complex $Y_i$,
a bonding map $\omega^{i+1}_i : Y_{i+1} \lo Y_i$ and
a map $\phi_i : X \lo Y_i$.
 We fix  metrics in $X$ and   in  each  $Y_i$
and with respect to these metrics we  determine  $0<\epsilon_i< 1/2^i$
 such  that
 the following properties will be satisfied:
 \\

 (i) $\phi_i$ is $0$-dimensional on $A$ and
 for every open  set $U \subset Y_i$ with $\diam U <2\epsilon_i$
 the set  $\phi_i^{-1}(U)\cap A_i $ splits into disjoint sets open in $A_i$ and
 of $\diam \leq 1/i$;

 (ii) $\dist(\omega^{i+1}_j\circ \phi_{i+1},\omega^{i}_j \circ \phi_{i}) < \epsilon_j/2^{i}$
 for $i\geq j$ where $\omega^j_i=\omega^j_{j-1} \circ \dots \circ \omega^{i+1}_i : Y_j \lo Y_i$
for $j>i$ and $\omega^i_i=id : Y_i \lo Y_i$.
 \\
 \\
 The construction will be carried out so that for $Y =\invlim (Y_i,\omega^{i+1}_i)$
  we  have  $\dim_\q Y \leq m-1$.
 Let us first show that the proposition follows from this construction.
 Denote  $f_i=\lim_{j\rightarrow \infty}  \omega^j_i\circ\phi_j : X \lo Y_i$.
 From (ii) it follows  that $f_i$ is well-defined, continuous and
 $\dist(f_i, \phi_i) \leq \epsilon_i$.
 From the definition of $f_i$ it follows that $f_j \circ f^j_i =f_i$.
 Hence the maps $f_i$ define the corresponding  map $f : X \lo Y$ such that
 $\omega_i \circ f=f_i$  where $\omega_i : Y \lo Y_i$ is the projection.
 Then it follows from (i) that for every $y \in Y_i$ the set
 $f_i^{-1}(y)\cap A_i$ splits into  finitely many  disjoint
 sets closed  in $A_i$ and of $\diam \leq 1/i$. This implies
 that  for every $y \in Y$  we have that $\dim( f^{-1}(y) \cap A_i)\leq 0$
 and hence $\dim (f^{-1}(y) \cap A)  \leq 0$.
 Then, by Corollary \ref{union-decomposition-for-dimension-types},
 $d_f \leq D_2 +1$. Since $\dim Y_i \leq m$ we have $\dim Y \leq m$ and
 since $\dim D_2 +1 =n-m$,  Hurewicz Theorem implies that $\dim Y = m$.
The condition  $\dim_\q Y \leq m-1$ and the formula for the cohomological dimension with respect
to $\z_{(p)}$
 imply that $\dim_{\z_{p\infty}}Y\le m-1$ for both $p$-regular and $p$-singular cases.
Therefore, $d_Y\leq D_1$. Thus, for every
$y \in Y$ we have
$$\dim (f^{-1}(y)\times Y) \leq \dim (D_1\boxplus (D_2 +1)) =n-1$$ and
hence $f$ is an exotic map.
 \\
 \\
 Now we return to our construction. Let $Y_1$  be an $m$-simplex
 and let $\phi_1 : X \lo Y$ be any map which is $0$-dimensional on $A$.
 Then one can choose $0<\epsilon_1< 1/2 $ so that (i) holds for $i=1$. Assume that
 the construction is completed for  $i$ and  proceed to $i+1$ as follows.
 Take a   triangulation of $Y_i$ so fine that for every simplex $\Delta$ of $Y_i$
 we have that $\diam (\omega^i_j (\Delta)) < \epsilon_j/2^i$ for every $j$ such that
 $i\geq j\geq 1$.
 For each $m$-simplex $\Delta$ of $Y_i$ consider a small  ball $D$ centered
 at the barycenter  of $\Delta$ and not touching $\partial D$.
 Recall that $\dim_\q X =m-1$. Then,
 by Proposition \ref{rational-dimension},
 for $\phi_i|_{\dots} :\phi_i^{-1}(\partial D)\lo \partial D$ there is
 a map $\psi_{\partial D} : \partial D \lo \partial D$ of non-zero degree such that
 $\psi_{\partial D} \circ \phi_i|_{\dots}$ extends over $\phi_i^{-1}(D)$ to a map
 $\phi_D : \phi_i^{-1}(D) \lo \partial D$. Clearly we can assume that
 $\psi_{\partial D}$ is $0$-dimensional (even finite-to-one).
 Denote by $\tilde Y_{i+1}$
 the quotient  space  of $Y_i$  obtained  by removing for each $m$-simplex
 $\Delta$ of $Y_i$ the interior $\Int D$ of the ball $D$ and
 identifying the points of  $\partial D$ according to the map $\psi_{\partial D}$.
  We  consider $\partial \Delta$
also as a subset of $\tilde Y_{i+1}$ and we denote by $Y_\Delta$ the subspace of
$\tilde Y_{i+1}$
obtained from  $\Delta$.
 Let
 $\tilde\omega_i^{i+1} : \tilde Y_{i+1} \lo Y_i$ be any map such
 that $\tilde \omega_i^{i+1}$ is $1$-to-$1$ over
 the $(m-1)$-skeleton of $Y_i$ and
 $\tilde\omega_i^{i+1}$ sends $Y_\Delta$ to $\Delta$ for every $m$-simplex $\Delta$ of $Y_i$.
 The map $\phi_i$ and the maps $\phi_D$ naturally define the corresponding  map
 $\tilde \phi_{i+1} : X \lo \tilde Y_{i+1}$ which coincides
 with  $\phi_i$ on $\phi_i^{-1}(\partial \Delta)$ for every $m$-simplex $\Delta$ of $Y_i$.
 Note that $\tilde \phi_{i+1}$ is $0$-dimensional on $A\cap \phi_{i}^{-1}(\Delta \setminus \Int D)$
 for every $m$-simplex $\Delta $ in $Y_i$.
 Let  a space $Y_{i+1}\supset \tilde Y_{i+1}$ be
 obtained from $\tilde Y_{i+1}$ by attaching for every sphere
 $\psi_{\partial D} (\partial D)\subset \tilde Y_{i+1}$
  the manifold
 $\psi_{\partial D} (\partial D) \times [0,1]$ by identifying
 $ \psi_{\partial D} (\partial D) \times 0$ with $\psi_{\partial D} (\partial D)$
 and let $\pi_{i+1} : Y_{i+1} \lo \tilde Y_{i+1}$ be a retraction
 projecting each manifold $\psi_{\partial D} (\partial D) \times [0,1]$ to
 $\psi_{\partial D} (\partial D)$. Then, by Proposition \ref{known},
 for every $m$-simplex $\Delta$ of $Y_i$ we can extend $\tilde \phi_{i+1}$ restricted
 to $\phi^{-1}_i(\Delta \setminus \Int D)$ over $\phi^{-1}_i(\Delta)$
 to a map being $0$-dimensional on $A \cap \phi^{-1}_i(\Delta)$ and
 this way we define the map $\phi_{i+1} : X \lo Y_{i+1}$ which
 is $0$-dimensional on $A$.
 Then there is $0<\epsilon_{i+1} < 1/2^{i+1}$
 such that (i) holds for $i+1$.
 Define $\omega_i^{i+1}=\tilde\omega_i^{i+1}\circ\pi_i: Y_{i+1} \lo Y_i$
 and note that (ii) is satisfied.
   Clearly  we may assume that $Y_{i+1}$ admits a triangulation and
 the construction is completed.

Note that for every $m$-simplex $\Delta$ in $Y_i$ we have that
 $$H^m((\omega^{i+1}_i)^{-1}(\Delta), (\omega^{i+1}_i)^{-1}(\partial \Delta);\q)=$$
$$H^m((\tilde\omega^{i+1}_i)^{-1}(\Delta), (\tilde\omega^{i+1}_i)^{-1}(\partial \Delta);\q)=
 H^m(Y_\Delta, \partial \Delta; \q )=0.$$
 This implies
 that for every subcomplex $Z_i$ of  $Y_i$ we have
 that $$H^m(Y_{i+1}, (\omega^{i+1}_i)^{-1}(Z_i); \q)=0.$$ Recall that
 $\diam \omega^{i}_j(\Delta)< \epsilon_j/2^i<1/2^{j+i}$ for every $i\geq j \geq 1$
 and every simplex $\Delta$ of $Y_i$.
 Then  $H^m(Y, Z;\q)=0$ for every closed subset $Z$ of $Y$ and
 hence $\dim_\q Y \leq m-1$.
 The proposition
 is proved. $\black$

\end{section}

Alexander Dranishnikov\\
Department of Mathematics\\
University of Florida\\
444 Little Hall\\
Gainesville, FL 32611-8105\\
dranish@math.ufl.edu\\
\\
Michael Levin\\
Department of Mathematics\\
Ben Gurion University of the Negev\\
P.O.B. 653\\
Be'er Sheva 84105, ISRAEL  \\
 mlevine@math.bgu.ac.il\\\\

\begin{thebibliography}{99}
 \bibitem{Bolt}  Boltyanskii, V.
 {\em An example of a two-dimensional compactum
 whose topological square is three-dimensional.} (Russian)
Doklady Akad. Nauk SSSR (N.S.) 67, (1949). 597-599 (English translation in
Amer. Math. Soc. Translation 1951, (1951). no. 48, 3--6).

\bibitem{DrUsp} Dranishnikov, A. N. {\em Homological dimension theory.} Russian Math. Surveys 43 (4) (1988), 11-63.

\bibitem{DrExt} Dranishnikov, A. N. {\em An extension of mappings into CW complexes.} Math. USSR Sb. 74 (1993), 47-56.

\bibitem{DrPacific} Dranishnikov A. N. {\em On the mapping intersection problem.} Pacific J. Math. 173 (1996), 403-412.


   \bibitem{dranishnikov-intersection}
   Dranishnikov, A. N. {\em On the dimension of the product of two compacta
   and the dimension of their intersection in general position in Euclidean space.}
   Trans. Amer. Math. Soc.  352  (2000),  no. 12, 5599--5618.

 

\bibitem{DrArxiv} Dranishnikov, A. N. {\em Cohomological dimension theory of compact metric spaces,} Topology
Atlas invited contribution, http://at.yorku.ca/topology.taic.html (see also  arXiv:math/0501523).

\bibitem{DRS1}
Dranishnikov A. N., Repovs D., Shchepin E.V. {\em On approximation and embedding problems for cohomological dimension.}
Topology Appl. 55 (1994), 67-86.

\bibitem{DRS}
Dranishnikov A. N., Repovs D., Shchepin E.V. {\em Transversal intersection formula for compacta.}
Topology Appl. 85 (1998), 93-117.

\bibitem{Dy}
 Dydak, Jerzy {\em Cohomological dimension and metrizable spaces. II.} Trans. Amer. Math. Soc. 348 (1996), no. 4, 1647--1661.

\bibitem{DyA}
Dydak, Jerzy {\em Algebra of dimension theory.} Trans. Amer. Math. Soc. 358 (2006), no. 4, 1537--1561.

\bibitem{Ko}
Kodama, Y. {\em Test spaces for homological dimension.}
Duke Math. J. 29 (1962) 41--50.



\bibitem{Ku}
Kuzminov, V. I. {\em Homological dimension theory.} Russian Math Surveys 23 (5)  (1968), 1-45.

   \bibitem{levin-lewis}
      Levin, Michael; Lewis, {\em Wayne Some mapping theorems for extensional dimension.}
      Israel J. Math.  133  (2003), 61-76.

   \bibitem{mill-pol}   van Mill, Jan; Pol, Roman {\em An example concerning
   the Menger-Urysohn formula.}
   Proc. Amer. Math. Soc.  138  (2010),  no. 10, 3749-3752.

\bibitem{Ol} Olszewski W. {\em Completion theorem for cohomological dimensions.} Proc. Amer. Math. Soc. 123 (1995) 2261--2264.

\bibitem{Pon} Pontryagin L. S. {\em Sur une hypothese fondamentale de la theorie de la dimension.}
C. R. Acad. Sci. Paris 190 (1930), 1105-1107.

\bibitem{Sch}
Shchepin, E. V. {\em Arithmetic of dimension theory.} Russian Math. Surveys 53 (1998), no. 5, 975--1069.


\end{thebibliography}
\end{document}